\documentclass[a4paper,10pt]{article}
\usepackage{amsmath,amssymb,graphicx}

\renewcommand{\title}[1]{\vspace{\fill}
\eject\addtolength{\baselineskip}{4pt}
{\bfseries\LARGE #1}\\[3mm]\addtolength{\baselineskip}{-4pt}}
\renewcommand{\author}[3]{\parbox[t]{75mm}
{\begin{center}{\scshape #1}\\[3mm] #2\\
 {\ttfamily #3} \end{center}}}

\newtheorem{thm}{\bfseries Theorem}
\newtheorem{lem}[thm]{\bfseries Lemma}        
\newtheorem{cor}[thm]{\bfseries Corollary}     

\newtheorem{cl}[thm]{\bfseries Claim}

\def\R{\mbox{\ensuremath{\mathbb R}}\xspace}

\def\F{\mbox{\ensuremath{\mathcal F}}\xspace}

\usepackage{xspace}
\usepackage{tikz}
\usetikzlibrary{calc}
\usepackage{hyperref}

\begin{document}

\begin{center}

\title{Improvement on line transversals of families of connected sets in the plane} 
\author{\underline{Miklós Csizmadia} \footnotemark[1] 
}{
Eötvös Loránd University \\
1117 Budapest, Pázmány Péter sétány 1/A, Hungary
}{
csizmadia.miki@gmail.com
} \footnotetext[1]{Supported by the EXCELLENCE-24 project no.~151504 Combinatorics and Geometry of the NRDI Fund.}


\end{center}


\begin{quote}
{\bfseries Abstract:}
Three lines are concurrent if they intersect at a single point. In this paper I prove that if \F is a bounded family of compact connected sets in the plane, such that every three sets in \F can be pierced by a single line, then there exists three concurrent lines in the plane such that the union of the three lines intersect every member of \F. This had previously only been proven for lines that are not required to be concurrent by McGinnis and Zerbib in \cite{ref1}. In fact, I prove a more general, ``colorful'' version of this result: If $\F_1, \dots , \F_5$ are bounded families of compact connected sets in the plane such that every three sets, chosen from three distinct families $\F_i$, can be pierced by a single line, then there exists $1 \leq j \leq 5$ and three concurrent lines, such that the union of the three lines intersect every member of $\F_j$. McGinnis and Zerbib had 6 families instead of 5, so I also improve their result in this respect. Moreover, the result can also be extended to unbounded families, if we allow the piercing lines to be parallel.



\end{quote}

\begin{quote}
{\bf Keywords:  }
combinatorial geometry, line transversals, line-piercing number, convex sets
\end{quote}
\vspace{5mm}


\section{Introduction} 

Let \F be a family of sets in the plane. We say that \F has the $T(r)$ property, if for every $r$ or fewer sets in \F, there exists a line, which intersects them. In other words, we say that every $r$ or fewer sized subset of \F \textit{admit a line transversal}. We say that \F is \textit{pierced} by $k$ lines, if there exists $k$ lines whose union intersect all sets in \F. The \textit{line-piercing number} of \F is the lowest such $k$.\\
The line-piercing number of sets in the plane has been a frequently studied topic since the 1960s. In particular, bounding the line-piercing number of compact convex sets proved to be an interesting problem. In 1964 Eckhoff \cite{ref6} proved that if a family of compact convex sets has the $T(4)$ property, then it can be pierced by 2 lines. In 1974 \cite{ref5} he gave an example of sets satisfying the $T(3)$ property which cannot be pierced by 2 lines. \\
The sets investigated in these results are assumed to be convex, but the results also apply to connected sets.  This follows from the fact that if $F$ is a connected set in $\R^2$, then a line $\ell$ intersects $F$ if and only if $\ell$ intersects $conv(F)$.\\
For a while, it was unknown whether the $T(3)$ property imply a finite line-piercing number. However, in 1975 Kramer \cite{ref7} proved that compact convex sets satisfying the $T(3)$ property can be pierced by 5 lines. After 18 years, in 1993 Eckhoff \cite{ref4} was able to show that such families can be pierced by 4 lines. Finally, in 2021 McGinnis and Zerbib \cite{ref1} proved that they can be pierced by only 3 lines. By their result, the line-piercing number of such families was finally resolved:
\begin{thm}[McGinnis and Zerbib \cite{ref1}]\label{thm1}
Let \F be a family of compact connected sets in $\R^2$. If every three sets in \F admit a line transversal, then \F is pierced by 3 lines.
\end{thm}

In fact, McGinnis and Zerbib also proved a generalized, ``colorful'' version of this result. In their paper they proved the theorem with 6 colors. Let $[n]$ denote the set $\{1,2, \dots, n\}$.
\begin{thm}[McGinnis and Zerbib \cite{ref1}]\label{thm2}
Let $\F_1, \dots , \F_6$ be families of compact connected sets in $\R^2$. If every three sets $F_1\in\F_{i_1},  \; F_2\in\F_{i_2}, \;  F_3\in\F_{i_3}, \quad 1\leq i_1 < i_2 < i_3 \leq 6$  admit a line transversal, then there exists $i\in [6]$ such that $\F_i$ can be pierced by 3 lines.
\end{thm}
This theorem generalizes Theorem \ref{thm1}, because if we choose all families to be the same family \F, it yields Theorem \ref{thm1}.\\
It can also be easily seen that if the statement is true for $k$ colors (families), then it is also true for all $l\geq k$ colors, by applying the $k$-color theorem to the first $k$ families. Therefore, lowering the number of colors makes a stronger statement.\\


\section{Main result} 

Here I will prove that Theorem \ref{thm2} is also true for concurrent lines and with 5 families, if the families are bounded. This theorem then can be extended to unbounded families if we allow the piercing lines to be parallel. Overall therefore, I also improve Theorem \ref{thm2} to 5 colors. \\
My proof is a modification of McGinnis and Zerbib's marvelous method for proving Theorem \ref{thm2}.

\begin{thm}\label{thm:main}
Let $\F_1, \dots , \F_5$ be bounded families of compact connected sets in $\R^2$. If every three sets $F_1\in\F_{i_1},  \; F_2\in\F_{i_2}, \;  F_3\in\F_{i_3}, \quad 1\leq i_1 < i_2 < i_3 \leq 5$  admit a line transversal, then there exists $i\in [5]$ such that $\F_i$ can be pierced by 3 concurrent lines.
\end{thm}
\begin{cor}
Let \F be a bounded family of compact connected sets in $\R^2$. If every three sets in \F admit a line transversal, then \F can be pierced by 3 concurrent lines.
\end{cor}

The main tool used in the proof will be the colorful KKM theorem (Knaster–Kuratowski–Mazurkiewicz lemma \cite{ref2}, generalized by David Gale \cite{ref3}).

\begin{thm}[KKM theorem \cite{ref2}]\label{kkm}
Let $\Delta^{n-1}$ be an $n-1$ dimensional simplex with $n$ vertices $\{v_1, \dots , v_n\}$. Let $A_1, \dots , A_n$ be open sets of $\Delta^{n-1}$ such that for every face $\sigma$ of $\Delta^{n-1}$, we have $\sigma \subset \bigcup^{}_{v_i\in \sigma}A_i$. Then $\bigcap A_i \neq \emptyset$. Namely, there is a point $x\in\Delta^{n-1}$ colored by all colors.
\end{thm}
In this case, we say that $A_1, \dots , A_n$ form a \textit{KKM-covering} of $\Delta^{n-1}$.

\begin{thm}[colorful KKM theorem \cite{ref3}]\label{ckkm}
Let $\Delta^{n-1}$ be an $n-1$ dimensional simplex. Let $A_1^i, \dots , A_n^i$, $i\in[n]$ be open sets of $\Delta^{n-1}$ such that for every $i\in[n]$, $\; A_1^i, \dots , A_n^i$ form a KKM-covering of $\Delta^{n-1}$. Then there exists $\pi\in S_n$ permutation, such that $\bigcap^{n}_{i=1}A_i^{\pi(i)} \neq \emptyset$.
\end{thm}

\medskip

{\medskip                    
\noindent{\scshape Proof of Theorem \ref{thm:main}:}}\\

First, we will prove a simplified version of this result, where the sets have a non-empty interior of size $\delta$ for some $\delta >0$. Formally, there exists $\delta>0$ such that each $F\in\F_i$ contain some translation of $B_{\delta}$ as a subset, where $B_{\delta}$ is the disk with radius $\delta$ centered at the origin. This condition applies in the special case where we have finite families of sets with non-empty interior, along with the condition of boundedness. Later, we will extend this proof to allow all compact connected sets.\\
We will also indirectly assume that no family can be pierced by 3 lines, and will arrive at a contradiction.\\
Using that the families are bounded, we may scale the plane such that every set in every family is contained inside the unit circle $U$. To use the KKM theorem, we will associate every point of a simplex with a set of 3 lines in the plane. Let $\Delta^4$ be the 4 dimensional simplex defined as follows: $\Delta^4 := \{(x_1, \dots ,x_5)\in\R^5 \;\; |\;\; x_i\geq 0,\ \sum x_i=1\}$ which is the convex hull of the canonical basis vectors $e_1, \dots , e_5$ of $\R^5$. We also parametrize the unit circle as $f(t) := (\cos(2\pi t), \sin(2\pi t))$. With this parametrization, a point  $x\in\Delta^4$ will correspond to 5 points $f_1(x), \dots , f_5(x)$ on $U$ in the following way: $f_i(x) := f(\sum_{j=1}^{i} x_j)$. Note that $f_5(x) = (1,0)$. \vspace{2pt}

Let $\overline{AB}$ denote the line defined by points $A$ and $B$. \vspace{2pt}
With the 5 points on the unit circle, we can define 3 lines: $\ell_1 := \overline{f_1(x)f_4(x)}, \;\; \ell_2 := \overline{f_2(x)f_5(x)}$, \; and \; $\ell_3 := \overline{f_3(x)M}$ where $M$ is the intersection of $\ell_1$ and $\ell_2$. The definition of $\ell_3$ is the main difference between this proof and McGinnis and Zerbib's. In their paper, they used a $6$th point on $U$ to define $\ell_3$.\\
Here however, $\ell_3$ is not defined when $\ell_1$ and $\ell_2$ coincide. We can avoid this, and other edge cases, by considering a slightly smaller simplex where each coordinate must be at least $\varepsilon$. Specifically, let $\Delta_{\varepsilon}^4 := \{(x_1, \dots ,x_5)\in\R^5 \;\; |\;\; x_i\geq \epsilon,\ \sum x_i=1\}$. This way, if $i \neq j$, then $f_i(x) \neq f_j(x)$. From now on, we will only consider $\Delta_{\varepsilon}^4$. The value of $\varepsilon$ will be determined later.\\

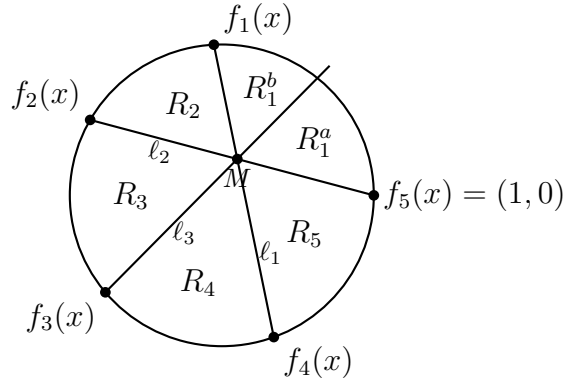
\begin{figure}[ht]
\begin{center}
\begin{tikzpicture}
    
    \def\r{2}
    
    \draw[thick] (0,0) circle (\r);
    
    \coordinate (P1) at (93:\r);
    \coordinate (P2) at (150:\r);
    \coordinate (P3) at (220:\r);
    \coordinate (P4) at (290:\r);
    \coordinate (P5) at (360:\r);

    \fill (P1) circle (2pt);
    \fill (P2) circle (2pt);
    \fill (P3) circle (2pt);
    \fill (P4) circle (2pt);
    \fill (P5) circle (2pt);

    \draw[thick] (P1) -- (P4);
    \draw[thick] (P2) -- (P5);

    \coordinate (M) at (intersection of P1--P4 and P2--P5);
    \fill (M) circle (2pt); 
    \node[anchor=north] at (M) {$M$};

    \draw[thick] (P3) -- ($ (P3)!1.7!(M) $);
    
    \node[anchor=south west] at (P1) {\large $f_1(x)$};
    \node[anchor=south east] at (P2) {\large $f_2(x)$};
    \node[anchor=north east] at (P3) {\large $f_3(x)$};
    \node[anchor=north west] at (P4) {\large $f_4(x)$};
    \node[anchor=west] at (P5) {\large $f_5(x) = (1,0)$};

    \node at (30:1.4) {\large $R_1^a$};
    \node at (70:1.5) {\large $R_1^b$};
    \node at (113:1.3) {\large $R_2$};
    \node at (180:1.2) {\large $R_3$};
    \node at (255:1.2) {\large $R_4$};
    \node at (335:1.2) {\large $R_5$};

    \node at (225:0.7) {$\ell_3$};
    \node at (310:1) {$\ell_1$};
    \node at (145:1) {$\ell_2$};

\end{tikzpicture}
\begin{minipage}{0.83\textwidth}
\caption{
A point $x\in\Delta_{\varepsilon}^4$ defines 5 points on the Unit circle, the 3 defined lines of which splits the unit disk into 6 regions.
}
\end{minipage}
\label{fig:first}
\end{center}
\end{figure}



Name each open region bounded by the lines and $U$ as $R_1^a, R_1^b, R_2, R_3, R_4, R_5$ as seen in Figure 1. Of course, the regions depend on $x$.\\
Now for every $1\leq j \leq 5$ we will define the colors corresponding to the family $\F_j$. For $i = 2,3,4,5$,  let $A_i^j$ be the set of $x\in\Delta_{\varepsilon}^4$ points, for which exists $F\in\F_j$ such that $F\subset R_i$.\: The case of $i=1$ is special: let $A_1^j$ be the set of $x\in\Delta_{\varepsilon}^4$ points, for which exists $F\in\F_j$ such that $F$ is the subset of either $R_1^a$ or $R_1^b$, this means $F$ is in one of the regions divided by $\ell_3$.
\begin{cl}
For each j, $A_1^j, \dots , A_5^j$ form a KKM cover of $\Delta_{\varepsilon}^4$.
\end{cl}
\textit{Openness:} \\
Since each $F\in \F_j$ is closed and each region is open, if $R_i$ contains a set $F$, then for a small ball around $x$, $R_i$ will still contain $F$, therefore all $A_i^j$ are open sets of $\Delta_{\varepsilon}^4$. Note that the regions $R_i$ change with $x$, and the set $F$ remains the in same position.\\
\textit{Condition on the faces of $\Delta_{\varepsilon}^4$:} \\
We will prove this by first showing that each face only has colors which are associated with the vertices of the face, then we will show that every point on the simplex is colored. \\
Firstly we need to check that if the $i$th coordinate of $x$ is $\varepsilon$, that is, we are on a face of $\Delta_{\varepsilon}^4$ not containing the $i$th vertex, then it cannot be colored by the $i$th color, formally $(x_i=\varepsilon) \; \Rightarrow \; (x \notin A_i^j)$. To prove this, we need to set $\varepsilon$ small enough that if $x_i=\varepsilon$, then $R_i$ cannot contain any set $F\in\F_j$. Here, we will use the fact that there exists $\delta>0$ such that each $F\in \F_j$ contains some translation of $B_{\delta}$ as a subset. Let $\varepsilon := \frac{\delta}{4\pi}$. With this $\varepsilon$, if $x_i = \varepsilon$, then a simple calculation shows that the length of the arc bounding $R_i$ must be less than $\delta$. Knowing this, using the fact that the bounding lines meet at $M$ inside $U$, we can conclude that $R_i$ cannot contain any set, therefore $(x_i=\varepsilon) \; \Rightarrow \; (x \notin A_i^j)$ holds.\\
Secondly, we have to show that, $\Delta_{\varepsilon}^4 = \bigcup_{i=1}^5 A_i^j$ for each $j$. This is true due to the indirect assumption that none of the families can be pierced by 3 lines. Suppose there was $x\in \Delta_{\varepsilon}^4$ such that $x\notin \bigcup_{i=1}^5 A_i^j$, then there is no $F\in\F_j$ which is in one of the 6 regions. Because every $F\in\F_j$ is connected, this means the union of $\ell_1, \ell_2, \ell_3$ pierce every $F\in\F_j$ which contradicts the assumption.
Knowing $(x_i=\varepsilon) \; \Rightarrow \; (x \notin A_i^j)$ and $\Delta_{\varepsilon}^4 = \bigcup_{i=1}^5 A_i^j$ we can conclude that the conditions of Theorem \ref{ckkm} hold. \hfill $\Box$\\
Thus, by Theorem \ref{ckkm}, there exists some permutation $\pi \in S_5$ and a point $p\in\Delta_{\varepsilon}^4$ such that $p\in\bigcap_{i=1}^{5} A_i^{\pi(j)}$. Therefore, for point $p$, each of the open regions $R_i$ contain a set $F_i \in \F_{\pi(i)}$, $i=2,3,4,5$. For $i=1$ it means either $R_1^a$ or $R_1^b$ contain a set $F_1\in\F_{\pi(1)}$.
If $F_1\subset R_1^a$, then the sets $F_1, F_2, F_4$ do not admit a line transversal. If $F_1\subset R_1^b$, then the sets $F_1, F_3, F_5$ do not admit a line transversal. In both cases, we arrive at a contradiction.
Therefore, the assumption that no $\F_j$ family can be pierced by 3 lines was incorrect, and there is a point $p\in\Delta_{\varepsilon}^4$ which is not colored by any $A_i^j$ for some $j$. The defined lines $\ell_1, \ell_2, \ell_3$ for $p$ then pierce the family $\F_j$. We can also see that these lines intersect at point $M$. \hfill $\Box$\\

With this, we can now generalize the proof for families of any compact connected sets:\\
Let $\F(\delta)$ be the thickened version of $\F$ by $\delta$. More formally, 
$\F(\delta) := \{ F + B_{\delta}\ | F\in\F \}$ for some $\delta \geq 0$, where $B_{\delta}$ is the closed disk with radius $\delta$ centered at the origin, and $F + B_{\delta}$ is the Minkowski sum of $F$ and $B_{\delta}$.
\begin{lem}\label{lem1}
Let \F be a bounded family of compact sets in the plane. If for all $\delta > 0$, $\F(\delta)$ is pierced by $n$ concurrent lines such that the intersection of said lines is always inside a fixed compact region, then \F is pierced by $n$ concurrent lines.
\end{lem}

We can see in the above proof, that the intersection $M$ is always inside the unit disk, therefore we can use Lemma \ref{lem1} in the context of 5 families to prove Theorem \ref{thm:main}. With the proof of this Lemma, the proof of Theorem \ref{thm:main} will be concluded.\\

\noindent{\scshape Proof of Lemma \ref{lem1}:}\\

Let us choose a sequence $(\delta_i)_{i=1}^\infty$, such that $\delta_i \to 0$ as $i\to\infty$. For every $\delta_i$, choose $n$ concurrent lines which pierce $\F(\delta_i)$: $\ell_1^i, \dots, \ell_n^i$. These lines intersect at point $M^i$. Let us also choose a point unit distance away from $M^i$ on each line: $P_1^i, \dots, P_n^i$. Using the fact that all $M^i$ and all $P_1^i, \dots, P_n^i$ are inside a fixed compact region, there must be a subsequence $(\delta_{i_k})_{k=1}^\infty$ $(i_1<i_2<i_3< \dots)$ such that all $M^{i_k}$ and all $P_1^{i_k}, \dots, P_n^{i_k}$ converge. Let these limit points be $M$ and $P_1, \dots, P_n$. These points define the lines $\ell_1, \dots, \ell_n$. We claim that $\F$ is pierced by $\ell_1, \dots, \ell_n$. Suppose, for contradiction, that there exists a set $F$ that is not pierced. Since $F$ is closed, there is a small enough $\delta_{i_K}$ such that $F(\delta_{i_K})$ is still not pierced. Since the lines $\ell_1^i, \dots, \ell_n^i$ converge to $\ell_1, \dots, \ell_n$, there is a big enough index $i_L>i_K$ such that $\ell_1^{i_L}, \dots, \ell_n^{i_L}$ still does not pierce $F(\delta_{i_K})$. Then, $\ell_1^{i_L}, \dots, \ell_n^{i_L}$ would not pierce $F(\delta_{i_L})$, which is a contradiction. \hfill $\Box$\\

With this, the proof of Theorem \ref{thm:main} is concluded. \hfill $\Box$\\

We can now extend Theorem \ref{thm:main} to unbounded families of compact connected sets:\\
\begin{lem}\label{lem2}
Let \F be a family of compact sets in the plane. If every bounded subfamily of \F is pierced by $n$ concurrent lines, then \F is pierced by $n$ concurrent or parallel lines.\\
\end{lem}

\noindent{\scshape Proof of Lemma \ref{lem2}:}\\

The proof is by induction on the number of lines.\\
Base case: $n=0$. This is obviously true.\\
Incuctive step: Assume the result holds for $n-1$, next we look at $n$.
If every bounded subfamily of $\F$ is pierced by $n-1$ concurrent lines, then by the inductive hypothesis we know that $\F$ can be pierced by $n-1$ concurrent lines. Hence, we can assume that there is a bounded subfamily which cannot be pierced by $n-1$ concurrent lines. Let us call this subfamily $\F_0$. Let $D_0$ be a large disk which amply contains $\F_0$. Let $D_i$ be the disk centered at the same point as $D_0$ but with radius increased by $i$ units ($i = 1,2,3,\dots$). Let $\F_i$ be the subfamily of $\F$ which is covered by $D_i$. Note that $\F_i \subseteq \F_j$ if $i \leq j$. We know that $\F_i$ can be pierced by $n$ concurrent lines: $\ell_1^i, \dots, \ell_n^i$ with intersection point $M^i$. Since $\F_0$ cannot be pierced by $n-1$ concurrent lines, each of $\ell_1^i, \dots, \ell_n^i$ must intersect $\F_0$, hence, each line must pass through $D_0$ as well. Let us choose two points one unit apart on each line, such that these points are inside $D_0$: $P_{1,1}^i, P_{1,2}^i, \dots, P_{n,2}^i$. This is possible because $D_0$ is sufficiently large. Because these points are inside $D_0$, there must be a subsequence $(i_k)_{k=1}^\infty$ along which all $2n$ points converge. Let the limit points be  $P_{1,1}, P_{1,2}, \dots, P_{n,2}$. These points define the limit lines: $\ell_1, \dots, \ell_n$. These lines will pierce every set in $\F$: each $F\in\F$ is pierced after some index $i_{K_F}$, and since these are closed sets, the limit lines will also pierce $F$. These limit lines will also be concurrent or parallel. This is true because if any two limit lines intersect at some point $M$, then the intersection point $M^{i_k}$ of the converging lines must converge to $M$, therefore in this case the limit lines are concurrent. The other case is when the limit lines are parallel, which can indeed happen if $M^{i_k}$ diverges unboundedly. With this, the lemma is now proven. \hfill $\Box$\\



Using Lemma \ref{lem2} the extension can now be formulated as follows:

\begin{thm}\label{}
Let $\F_1, \dots , \F_5$ be families of compact connected sets in $\R^2$. If every three sets $F_1\in\F_{i_1},  \; F_2\in\F_{i_2}, \;  F_3\in\F_{i_3}, \quad 1\leq i_1 < i_2 < i_3 \leq 5$  admit a line transversal, then there exists $i\in [5]$ such that $\F_i$ can be pierced by 3 concurrent or parallel lines.
\end{thm}
\begin{cor}
Let \F be a family of compact connected sets in $\R^2$. If every three sets in \F admit a line transversal, then \F can be pierced by 3 concurrent or parallel lines.
\end{cor}

Similarly to McGinnis and Zerbib's paper \cite{ref1}, the proof of Theorem \ref{thm:main} imply a slightly stronger result: if the families are bounded, due to the fact that $f_5(x)$ is fixed every time, we can choose the position of $f_5(x)$. Therefore, we can choose any point $Q$ for which one of the piercing lines will go through, as long as $Q$ is outside the convex hull of the sets.\\
Moreover, by choosing $Q$ further and further in one direction, a simple convergence argument shows, that the direction of one of the piercing lines can be chosen in advance.\\
All the results in this paper can be demonstrated in Eckhoff's \cite{ref5} example of sets which cannot be pierced by 2 lines. His example can indeed be pierced by 3 concurrent lines, one of whose direction can be chosen.


\subsubsection*{Acknowledgment:}
I am deeply grateful to my undergraduate supervisor, Dömötör Pálvölgyi, for his insightful comments on this paper and his guidance in my research.


\begin{thebibliography}{99}


\bibitem{ref1}
D. McGinnis, S. Zerbib, Line transversals in families of connected sets in the plane, SIAM Journal on Discrete Mathematics
36(4) (2022), 2916-2919. 


\bibitem{ref2}
B. Knaster, C. Kuratowski, and S. Mazurkiewicz, Ein Beweis des Fixpunktsatzes für n-Dimensionale Simplexe, Fund. Math., 14 (1929), 132-137.


\bibitem{ref3}
D. Gale and A. Mas-Colell, An equilibrium existence theorem for a general model without ordered preferences, Journal of Mathematical Economics, 2(1) (1975), 9–15.


\bibitem{ref4}
J. Eckhoff, A Gallai-type transversal problem in the plane. Discrete Comput. Geom. 9 (1993), no. 2, 203–214.


\bibitem{ref5}
J. Eckhoff, Transversatenprobleme in der Ebene, Arch. Math. 24 (1973), 195–202.


\bibitem{ref6}
J. Eckhoff, Transversalenprobleme vom Gallaischen Typ, Dissertation, Universität Göttingen,
(1969).

\bibitem{ref7}
D. Kramer, Transversalenprobleme vom Hellyschen und Gallaischen Typ, Dissertation, Universit¨at Dortmund, (1974).


\end{thebibliography}
\end{document}